\documentclass[12pt]{article}
\usepackage{amssymb,amsmath,latexsym,amsthm}
\usepackage{url}
\usepackage{enumerate}
\usepackage{ifpdf}
\ifpdf
  \usepackage[pdftex]{hyperref}
\else
  \usepackage[dvips]{hyperref}
\fi

\title{The calculus of differentials for the weak Stratonovich integral}
\author{Jason Swanson\thanks{Supported in part by NSA grant H98230-09-1-0079.}\\
University of Central Florida\\
\url{http://math.swansonsite.com/}
}
\date{July 31, 2011}

\usepackage{color}

\begin{document}

\newtheorem{thm}{Theorem}[section]
\newtheorem{cor}[thm]{Corollary}
\newtheorem{lemma}[thm]{Lemma}
\theoremstyle{definition}
\newtheorem{defn}[thm]{Definition}
\newtheorem{rmk}[thm]{Remark}
\newtheorem{expl}[thm]{Example}
\numberwithin{equation}{section}
\def\al{\alpha}
\def\be{\beta}
\def\Ga{\Gamma}
\def\de{\delta}
\def\De{\Delta}
\def\ep{\varepsilon}
\def\th{\theta}
\def\Th{\Theta}
\def\ka{\kappa}
\def\la{\lambda}
\def\La{\Lambda}
\def\si{\sigma}
\def\ph{\varphi}
\def\dd{\mathbf{d}}
\def\FF{\mathcal{F}}
\def\GG{\mathcal{G}}
\def\II{\mathcal{I}}
\def\MM{\mathcal{M}}
\def\RR{\mathbb{R}}
\def\SS{\mathcal{S}}
\def\VV{\mathcal{V}}
\def\ZZ{\mathbb{Z}}
\def\bu{\bullet}
\def\wt{\widetilde}
\def\wh{\widehat}
\def\ol{\overline}
\def\To{\Rightarrow}
\def\pf{\noindent{\bf Proof.} }
\def\qed{\hfill $\Box$}
\providecommand{\flr}[1]{\lfloor#1\rfloor}
\providecommand{\cub}[1]{[\![#1]\!]}
\providecommand{\ang}[1]{\langle#1\rangle}

\maketitle

\begin{abstract}

The weak Stratonovich integral is defined as the limit, in law, of Stratonovich-type symmetric Riemann sums. We derive an explicit expression for the weak Stratonovich integral of $f(B)$ with respect to $g(B)$, where $B$ is a fractional Brownian motion with Hurst parameter 1/6, and $f$ and $g$ are smooth functions. We use this expression to derive an It\^o-type formula for this integral. As in the case where $g$ is the identity, the It\^o-type formula has a correction term which is a classical It\^o integral, and which is related to the so-called signed cubic variation of $g(B)$. Finally, we derive a surprising formula for calculating with differentials. We show that if $\dd M=X\,\dd N$, then $Z\,\dd M$ can be written as $ZX\,\dd N$ minus a stochastic correction term which is again related to the signed cubic variation.

\bigskip


\bigskip


\end{abstract}

\section{Introduction}

If $X$ and $Y$ are stochastic processes, then the Stratonovich integral of $X$ with respect to $Y$ can be defined as the ucp (uniformly on compacts in probability) limit, if it exists, of the process
  \[
  t \mapsto \sum_{t_j\le t}\frac{X(t_{j-1}) + X(t_j)}2
    (Y(t_j) - Y(t_{j-1})),
  \]
as the mesh of the partition $\{t_j\}$ goes to zero. If we specialize to the uniformly spaced partition, $t_j=j/n$, then we are interested in the Stratonovich-type symmetric Riemann sums,
  \begin{equation}\label{Riem_sum}
  \sum_{j=1}^{\flr{nt}}
    \frac{X(t_{j-1}) + X(t_j)}2(Y(t_j)-Y(t_{j-1})),
  \end{equation}
where $\flr{x}$ denotes the greatest integer less than or equal to $x$.

It is well-known (see \cite{CN} and \cite{GNRV}) that if $Y=B^H$, a fractional Brownian motion with Hurst parameter $H$, and $X=f(B^H)$ for a sufficiently differentiable function $f$, then the Stratonovich integral of $X$ with respect to $Y$ exists for all $H>1/6$, but does not exist for $H=1/6$. Moreover, if $H>1/6$, then the Stratonovich integral satisfies the classical Stratonovich change-of-variable formula, which corresponds to the usual fundamental theorem of calculus.

In \cite{NRS}, we studied the case $H=1/6$. There we showed that if $Y=B=B^{1/6}$ and $X=f(B)$, where $f\in C^\infty(\RR)$, then the sequence of processes \eqref{Riem_sum} converges in law. We let $\int_0^t f(B(s)) \,\dd\mathbb{d} B(s)$ denote a process with this limiting law, and we referred to this as the weak Stratonovich integral. We also showed that the weak Stratonovich integral with respect to $B$ does not satisfy the classical Stratonovich change-of-variable formula. Rather, it satisfies an It\^o-type formula with a correction term that is a classical It\^o integral. Namely,
  \begin{equation}\label{ito_form}
  f(B(t)) = f(B(0)) + \int_0^t f'(B(s))\,\dd B(s)
    - \frac1{12}\int_0^t f'''(B(s))\,\dd\cub{B}_s,
  \end{equation}
where $\cub{B}$ is what we called the signed cubic variation of $B$. That is, $\cub{B}$ is the limit in law of the sequence of processes $\sum_{j=1}^{\flr{nt}}(B(t_j)-B(t_{j-1}))^3$. It is shown in \cite{NO} that $\cub{B}=\ka W$, where $W$ is a standard Brownian motion, independent of $B$, and $\ka$ is an explicitly defined constant whose approximate numerical value is $\ka\simeq 2.322$. (See \eqref{kappa} for the precise definition of $\ka$.) The correction term above is a standard It\^o integral with respect to Brownian motion. Similar It\^o-type formulas with an It\^o integral correction term were developed in \cite{BS} and \cite{NR}. There, the focus was on quartic variation processes and midpoint-style Riemann sums. A formula similar to \eqref{ito_form}, but with an ordinary integral correction term, was established in \cite{ER} for the Russo-Vallois symmetric integral with respect to finite cubic variation processes.

The precise results in \cite{BS} and \cite{NRS}, as well as in this paper, involve demonstrating the joint convergence of all of the processes involved, with the type of convergence being weak convergence as processes in the Skorohod space of c\`adl\`ag functions. In Section \ref{S:setup}, we establish the formal definition of the weak Stratonovich integral as an equivalence class of sequences of c\`adl\`ag step functions, and we demonstrate in Theorem \ref{T:gen_main} the joint convergence in law of such sequences. For simplicity, we omit discussion of these details in this introduction, and only summarize the results of Section \ref{S:formulas}, in which we derive our various change-of-variable formulas.

In Section \ref{S:formulas}, we extend the It\^o-type formula \eqref{ito_form} to the case $Y=g(B)$. We show that the sequence of processes \eqref{Riem_sum} converges in law to an integral satisfying the It\^o-type formula
  \begin{equation}\label{old_ito_form_Y}
  \ph(Y(t)) = \ph(Y(0)) + \int_0^t \ph'(Y(s))\,\dd Y(s)
    - \frac1{12}\int_0^t \ph'''(Y(s))\,\dd\cub{Y}_s,
  \end{equation}
where
  \[
  \cub{Y}_t = \int_0^t (g'(B(s)))^3\,\dd\cub{B}_s
  \]
is the limit, in law, of $\sum_{j=1}^{\flr{nt}}(Y(t_j)-Y(t_{j-1}))^3$. That is, $\cub{Y}$ is the signed cubic variation of $Y$.

This result is actually just one of the two main corollaries of our central result. (See Corollary \ref{C:ito_form_Y}.) To motivate the other results, consider the following. Formulas such as \eqref{ito_form} and \eqref{old_ito_form_Y} are typically referred to as change-of-variable formulas. They have the same structure as It\^o's rule, which is also generally referred to as a change-of-variable formula. In elementary calculus, we perform a change-of-variable when we convert an integral with respect to one variable into an integral with respect to another. In It\^o's stochastic calculus, we may wish to convert an integral with respect to one semimartingale into an integral with respect to another. Strictly speaking, It\^o's rule is not sufficient for this purpose. It\^o's rule simply tells us how to expand a function of a semimartingale into a sum of integrals. In order to convert one integral into another, we must combine It\^o's rule with a theorem that says:
  \[
  \text{If $M=\int X\,dY$, then $\int Z\,dM = \int ZX\,dY$.}
  \]
Or, in differential form:
  \begin{equation}\label{sub}
  \text{If $dM=X\,dY$, then $Z\,dM = ZX\,dY$.}
  \end{equation}
For It\^o integrals, this theorem is usually proved very early on in the construction of the integral. It is also true for the classical Stratonovich integral for semimartingales, as well as for ordinary Lebesgue-Stieltjes integrals. In fact, in the theory of Lebesgue-Stieltjes integration, it is often this result which is called the change-of-variable formula.

In terms of the calculus of differentials, It\^o's rule tells us that if $M=f(Y)$, then $dM=f'(Y)\,dY + \frac12f''(Y)\,d\ang{Y}$, where $\ang{Y}$ is the quadratic variation of $Y$; and \eqref{sub} tells us that it is permissible to substitute this expression into $Z\,dM$, so that $Z\,dM=Zf'(Y)\,dY + \frac12Zf''(Y)\,d\ang{Y}$.

In this paper, we will show that \eqref{sub} is not true for the weak Stratonovich integral. A very simple example which illustrates this is the following. First, let us note that when the integral is defined as a limit of Stratonovich-type symmetric Riemann sums, it is always the case that $\int\th\,d\th=\frac12\th^2$, for any process $\th$. Let us therefore define $M=\frac12B^2$, so that $\dd M =B\,\dd B$. On the other hand,
  \[
  \int M\,\dd M = \frac12 M^2 = \frac18 B^4.
  \]
Using \eqref{ito_form}, we have
  \[
  \frac18 B^4 = \int \frac12 B^3\,\dd B
    - \frac1{12}\int 3B\,\dd\cub{B}
    = \int MB\,\dd B - \frac14\int B\,\dd\cub{B}.
  \]
It follows that, in this example, \eqref{sub} does not hold for the weak Stratonovich integral. Instead, we have that $\dd M=B\,\dd B$, whereas $M\,\dd M = MB\,\dd B-\frac14B\,\dd\cub{B}$.

The second main corollary of our central result is that the weak Stratonovich integral satisfies a rule analogous to \eqref{sub}, but with a correction term. (See Corollary \ref{C:sub_weak_semi}.) Namely, suppose $X=f(B)$, $Y=g(B)$, and $Z=h(B)$, where $f,g,h\in C^\infty(\RR)$. Then the weak Stratonovich integral satisfies the following rule for calculating with differentials:
  \begin{equation}\label{sub_weak}
  \text{If $\dd M=X\,\dd Y$, then
    $Z\,\dd M = ZX\,\dd Y-\frac14(f'g'h')(B)\,\dd \cub{B}$.}
  \end{equation}
We actually prove a slightly more general rule; see \eqref{sub_weak_semi}.

Both \eqref{old_ito_form_Y} and \eqref{sub_weak} will be demonstrated as corollaries of the following general result. With $X$ and $Y$ as above,
  \begin{equation}\label{int_def}
  \int_0^t X(s)\,\dd Y(s) = \Phi(B(t)) - \Phi(B(0))
    + \frac1{12}\int_0^t(f''g' - f'g'')(B(s))\,\dd\cub{B}_s,
  \end{equation}
where $\Phi\in C^\infty(\RR)$ is chosen to satisfy $\Phi'=fg'$. See Theorem \ref{T:main} for the precise statement. Theorem \ref{T:main} is actually formulated more generally, for integrators of the form $Y+V$, where $V=\int\th(B)\,\dd\cub{B}$. This generalization is necessary to make sense of $\int Z\,\dd M$ in \eqref{sub_weak}, since if $M=\int X\,\dd Y$, then according to \eqref{int_def}, $M$ is not a function of $B$, but is rather the sum of a function of $B$ and a process $V$ which is in an integral against $\cub{B}$.

\section{Notation and definitions}\label{S:setup}

\subsection{Basic notation}

Let $B=B^{1/6}$ be a fractional Brownian motion with Hurst parameter $H=1/6$. That is, $B$ is a centered Gaussian process, indexed by $t\ge0$, such that
  \[
  E[B(s)B(t)] = \frac12(t^{1/3} + s^{1/3} - |t - s|^{1/3}).
  \]
For compactness of notation, we will sometimes write $B_t$ instead of $B(t)$, and similarly for other processes. Given a positive integer $n$, let $t_j=t_{j,n}=j/n$. We shall frequently have occasion to deal with the quantity
  \[
  \be_j = \be_{j,n} = \frac{B(t_{j-1}) + B(t_j)}2.
  \]
Let $\De B_{j,n}=B(t_j)-B(t_{j-1})$ and $B^*(T)=\sup_{0\le t\le T} |B(t)|$.

Let $\ka>0$ be defined by
  \begin{equation}\label{kappa}
  \ka^2 = \frac34\sum_{r\in\ZZ}
    (|r + 1|^{1/3} + |r - 1|^{1/3} - 2|r|^{1/3})^3.
  \end{equation}
Let $D_{\RR^d}[0,\infty)$ denote the Skorohod space of c\`adl\`ag functions from $[0,\infty)$ to $\RR^d$. Throughout the paper, ``$\To$" will denote convergence in law. The phrase ``uniformly on compacts in probability" will be abbreviated ``ucp." If $X_n$ and $Y_n$ are c\`adl\`ag processes, we shall write $X_n\approx Y_n$ or $X_n(t)\approx Y_n(t)$ to mean that $X_n-Y_n\to0$ ucp.

\subsection{The space $[\SS]$}

Recall that for fixed $n$, we defined $t_k=k/n$. Let $\SS_n$ denote the vector space of stochastic processes $\{L(t): t\ge0\}$ of the form $L = \sum_{k=0}^\infty \la_k 1_{[t_k,t_{k+1})}$, where each $\la_k\in\FF^B_\infty$. Note that $\la_k=L(t_k)$. Given $L\in\SS_n$, let $\de_j(L) = L(t_j) - L(t_{j-1})$, for $j\ge 1$. Since $t\in[t_k, t_{k+1})$ if and only if $\flr{nt}=k$, we may write
  \[
  L(t) = L(0) + \sum_{j=1}^{\flr{nt}} \de_j(L).
  \]

\begin{defn}\label{D:script_S}
Let $\SS$ denote the vector space of sequences $\La=\{\La_n\}_{n=1}^\infty$ such that
  \begin{enumerate}[(i)]
  \item $\La_n\in\SS_n$,
  \item $\La_n(0)$ converges in probability, and
  \item there exist $\ph_1,\ph_3,\ph_5 \in C^\infty(\RR)$ such that
    \begin{equation}\label{script_S}
    \de_j(\La_n) = \ph_1(\be_j)\De B_{j,n}
      + \ph_3(\be_j)\De B_{j,n}^3
      + \ph_5(\be_j)\De B_{j,n}^5 + R_{j,n},
    \end{equation}
  where for each $T,K>0$, there exists a finite constant $C_{T,K}$ such that
  \[
  |R_{j,n}|1_{\{B^*(T)\le K\}} \le C_{T,K}|\De B_{j,n}|^7,
  \]
whenever $j/n\le T$.
  \end{enumerate}
\end{defn}

If $X=f(B)$, where $f\in C^\infty(\RR)$, then we define
  \[
  \La^X_n = \sum_{k=0}^\infty X(t_k) 1_{[t_k,t_{k+1})},
  \]
and $\La^X = \{\La^X_n\}_{n =1}^\infty$. Note that the map $X\mapsto \La^X$ is linear.

\begin{lemma}
If $X=f(B)$, where $f\in C^\infty(\RR)$, then $\La^X\in\SS$ and $\La^X_n \to X$ uniformly on compacts a.s.
\end{lemma}

\pf Since $X$ is continuous a.s., we have that $\La^X_n\to X$ uniformly on compacts a.s. Clearly, $\La^X_n\in\SS_n$ and $\La^X_n(0)=X(0)$ for all $n$, so that Definition \ref{D:script_S}(i) and (ii) hold. For $a,b\in\RR$, we use the Taylor expansion
  \[
  f(b) - f(a) = f'(x)(b - a) + \frac1{24}f'''(x)(b - a)^3
    + \frac1{5!2^4}f^{(5)}(x)(b - a)^5
    + h(a,b)(b - a)^7,
  \]
where $x=(a+b)/2$, and $|h(a,b)|\le M(a,b)=\sup_{x\in[a\wedge b,a\vee b]} |g^{(7)}(x)|$. For a derivation of this Taylor expansion, see the proof of Lemma 5.2 in \cite{NRS}.

Taking $a=B(t_{j-1})$ and $b=B(t_j)$ gives
  \begin{equation}\label{smooth_expan}
  \begin{split}
  \de_j(\La^X_n) &= f(B(t_j)) - f(B(t_{j-1}))\\
  &= f'(\be_j)\De B_{j,n}
    + \frac1{24}f'''(\be_j)\De B_{j,n}^3
    + \frac1{5!2^4}f^{(5)}(\be_j)\De B_{j,n}^5
    + R_{j,n},
  \end{split}
  \end{equation}
where $|R_{j,n}| \le M(B(t_{j-1}),B(t_j))|\De B_{j,n}|^7$. If $j/n\le T$ and $B^*(T)\le K$, then $B(t_{j-1}),B(t_j)\in[-K,K]$, which implies $M(B(t_{j-1}),B(t_j))\le \sup_{x\in[-K,K]}|g^{(7)}(x)|<\infty$, and this verifies Definition \ref{D:script_S}(iii) showing that $\La^X \in\SS$. \qed

\bigskip

We may now identify $X=f(B)$ with $\La^X\in\SS$, and will sometimes abuse notation by writing $X\in\SS$. In this way, we identify the space of smooth functions of $B$ with a space of sequences in such a way that each sequence converges a.s. to its corresponding process. What we see next is the every sequence in $\SS$ converges to a stochastic process, at least in law.

\begin{thm}\label{T:gen_main}
Let $\La^{(1)},\ldots,\La^{(m)}\in\SS$. For $1\le k\le m$, choose $\ph_{1,k}, \ph_{3,k},\ph_{5,k} \in C^\infty(\RR)$ satisfying \eqref{script_S} for $\La^{(k)}$ and let $\II^{(k)}(0)$ be the limit in probability of $\La^{(k)}_n(0)$ as $n\to\infty$. Let $\Phi_k\in C^\infty(\RR)$ satisfy $\Phi_k'= \ph_{1,k}$ and $\Phi_k(0)=0$. Let $W$ be a Brownian motion independent of $B$, and let $\ka>0$ be given by \eqref{kappa}. Define
  \[
  \II^{(k)}(t) = \II^{(k)}(0) + \Phi_k(B(t))
    + \ka\int_0^t \left({
    \ph_{3,k} - \frac1{24}\ph_{1,k}''
    }\right)(B(s))\,dW(s),
  \]
where this last integral is an It\^o integral. Then $(B,\La^{(1)}_n, \ldots,\La^{(m)}_n)\To (B,\II^{(1)},\ldots,\II^{(m)})$ in $D_{\RR^{m+1}} [0,\infty)$ as $n\to\infty$.
\end{thm}

\pf By Definition \ref{D:script_S}, we may write
  \[
  \La^{(k)}_n(t) = \La^{(k)}_n(0)
    + \sum_{j=1}^{\flr{nt}}\ph_{1,k}(\be_j)\De B_{j,n}
    + \sum_{j=1}^{\flr{nt}}\ph_{3,k}(\be_j)\De B_{j,n}^3
    + \sum_{j=1}^{\flr{nt}}\ph_{5,k}(\be_j)\De B_{j,n}^5
    + R_n(t),
  \]
where $R_n(t)=\sum_{j=1}^{\flr{nt}}R_{j,n}$. let $R_n^*(T) = \sup_{0 \le t\le T}|R_n(t)|\le\sum_{j=1}^{\flr{nT}}|R_{j,n}|$. Let $\ep>0$ and choose $K$ such that $P(B^*(T)>K)<\ep$. Then
  \[
  P(R^*_n(T) > \ep) \le P(B^*(T) > K)
    + P\bigg(C_{T,K}\sum_{j=1}^{\flr{nT}}
    |\De B_{j,n}|^7 > \ep\bigg).
  \]
Since $B$ has a nontrivial 6-variation (see Theorem 2.11 in \cite{NRS}), we have $\sum_{j=1}^{\flr{nT}} |\De B_{j,n}|^7\to0$ a.s. Hence, for $n$ sufficiently large, we have $P(R^*_n(T) > \ep)<2\ep$, which gives $R_n\to0$ ucp.

As in the proof of Theorem 2.13 in \cite{NRS}, we may assume without loss of generality that each $\ph_{i,k}$ has compact support. By Lemma 5.1 in \cite{NRS}, if $\ph\in C^1(\RR)$ has compact support, then $\sum_{j=1}^{\flr{nt}}\ph(\be_j)\De B_{j,n}^5\to0$ ucp. Thus,
  \[
  \La^{(k)}_n(t) \approx \II^{(k)}(0)
    + \sum_{j=1}^{\flr{nt}}\ph_{1,k}(\be_j)\De B_{j,n}
    + \sum_{j=1}^{\flr{nt}}\ph_{3,k}(\be_j)\De B_{j,n}^3.
  \]
Similarly, by \eqref{smooth_expan},
  \begin{align*}
  \Phi_k(B(t))
    &\approx \sum_{j=1}^{\flr{nt}}
    (\Phi_k(B(t_j)) - \Phi_k(B(t_{j-1})))\\
  &\approx \sum_{j=1}^{\flr{nt}}\ph_{1,k}(\be_j)\De B_{j,n}
    + \frac1{24}\sum_{j=1}^{\flr{nt}}
    \ph_{1,k}''(\be_j)\De B_{j,n}^3.
  \end{align*}
Therefore,
  \[
  \La^{(k)}_n(t) \approx \II^{(k)}(0) + \Phi_k(B(t))
    + \sum_{j=1}^{\flr{nt}}\psi_k(\be_j)\De B_{j,n}^3,
  \]
where $\psi_k=\ph_{3,k}-\frac1{24}\ph_{1,k}''$. Let $V_n(\psi,t) = \sum_{j=1}^{\flr{nt}}\psi(\be_j) \De B_{j,n}^3$ and $J_k(t) = \ka\int_0^t\psi_k(B(s))\,dW(s)$. By Lemma 5.2 and Theorem 2.13 in \cite{NRS}, we have $(B,V_n(\psi_1),\ldots,V_n(\psi_m)) \To (B,J_1, \ldots, J_m)$,
in $D_{\RR^{m+1}} [0,\infty)$ as $n\to\infty$, which implies $(B, \La^{(1)}_n, \ldots,\La^{(m)}_n)\To (B,\II^{(1)},\ldots,\II^{(m)})$. \qed

\bigskip

We now define an equivalence relation on $\SS$ by $\La\equiv\Th$ if and only if $\La_n-\Th_n\to0$ ucp.

\begin{lemma}
If $\La\in\SS$, then there exist unique functions $\ph_1,\ph_3$ which satisfy \eqref{script_S}. If we denote these unique functions by $\ph_{1,\La}$ and $\ph_{3,\La}$, then $\La\equiv\Th$ if and only if
  \begin{enumerate}[(i)]
  \item $\La_n(0)-\Th_n(0)\to0$ in probability, and
  \item $\ph_{1,\La}=\ph_{1,\Th}$ and $\ph_{3, \La}=\ph_{3,\Th}$.
  \end{enumerate}
\end{lemma}

\pf Let $\La\in\SS$. Let $\{\ph_1,\ph_3,\ph_5\}$ and $\{\wt\ph_1, \wt\ph_3,\wt\ph_5\}$ be two sets of functions, each of which satisfies \eqref{script_S}. Let $\II(0)$ be the limit in probability of $\La_n(0)$ as $n\to\infty$. Let $\Phi,\wt\Phi\in C^\infty(\RR)$ satisfy $\Phi'= \ph_1$, $\wt\Phi'= \wt\ph_1$, and $\Phi(0)=\wt\Phi(0)=0$. Then, by Theorem \ref{T:gen_main}, $\La_n$ converges in law in $D_\RR[0,\infty)$ to
  \begin{multline*}
  \II(t) = \II(0) + \Phi(B(t)) + \ka\int_0^t \left({
    \ph_3 - \frac1{24}\ph_1''}\right)(B(s))\,dW(s)\\
    = \II(0) + \wt\Phi(B(t)) + \ka\int_0^t \left({
    \wt\ph_3 - \frac1{24}\wt\ph_1''}\right)(B(s))\,dW(s).
  \end{multline*}
Hence, $E[\II(t)-\II(0)\mid\FF_\infty^B]=\Phi(B(t))=\wt\Phi(B(t))$ a.s. for all $t\ge0$, which implies $\Phi=\wt\Phi$, and hence, $\ph_1=\wt\ph_1$. It follows that
  \[
  \MM(t) = \int_0^t (\ph_3 - \wt\ph_3)(B(s))\,dW(s) = 0.
  \]
Hence, $E[\MM(t)^2 \mid \FF_\infty^B] = \int_0^t |(\ph_3 - \wt\ph_3) (B(s))|^2\,ds=0$ a.s. for all $t\ge0$, which implies $\ph_3=\wt\ph_3$. This shows that there exist unique functions $\ph_{1,\La},\ph_{3,\La}$ which satisfy \eqref{script_S}.

Let $\La,\Th\in\SS$ and define $\Ga=\La-\Th$. Note that $\La_n - \Th_n\to0$ ucp if and only if $\Ga_n\To0$ in $D_\RR[0,\infty)$.

First assume (i) and (ii) hold. Then $\Ga_n(0)\to0$ in probability, so by Theorem \ref{T:gen_main}, $\Ga_n$ converges in law in $D_\RR[0,\infty)$ to
  \[
  \Phi_\Ga(B(t)) + \ka\int_0^t \left({
    \ph_{3,\Ga} - \frac1{24}\ph_{1,\Ga}''}\right)(B(s))\,dW(s),
  \]
where $\Phi_\Ga'=\ph_{1,\Ga}$ and $\Phi_\Ga(0)=0$. But from \eqref{script_S}, we see that $\ph_{1,\Ga}=\ph_{1,\La}-\ph_{1,\Th}=0$ and $\ph_{3,\Ga} =\ph_{3,\La}-\ph_{3,\Th}=0$. Hence, $\Ga_n\To0$ and $\La\equiv\Th$.

Now assume $\La\equiv\Th$. Then $\Ga_n\to0$ ucp, so by Theorem \ref{T:gen_main}, for all $t\ge 0$,
  \[
  \II(t) = \II(0) + \Phi_\Ga(B(t)) + \ka\int_0^t \left({
    \ph_{3,\Ga} - \frac1{24}\ph_{1,\Ga}''
    }\right)(B(s))\,dW(s) = 0,
  \]
where $\II(0)$ is the limit in probability of $\La_n(0)-\Th_n(0)$ as $n\to\infty$, and $\Phi_\Ga'=\ph_{1,\Ga}$ with $\Phi_\Ga(0)=0$. Thus, $\II(0)=0$, which shows that (i) holds. And as above, we obtain $\ph_{1,\Ga}=\ph_{3,\Ga}=0$, which shows that (ii) holds. \qed

\bigskip

Let $[\La]$ denote the equivalence class of $\La$ under this relation, and let $[\SS]$ denote the set of equivalence classes. If $N=[\La] \in[\SS]$, then we define $\ph_{1,N}=\ph_{1,\La}$, $\ph_{3,N} = \ph_{3,\La}$, $\II_N(0)=\lim\La_n(0)$, and
  \begin{equation}\label{I_def}
  \II_N(t) = \II_N(0) + \Phi_N(B(t))
    + \ka\int_0^t \left({
    \ph_{3,N} - \frac1{24}\ph_{1,N}''}\right)(B(s))\,dW(s),
  \end{equation}
where $\Phi_N'=\ph_{1,N}$ and $\Phi_N(0)=0$. Notice that by Theorem \ref{T:gen_main}, if $N_1,\ldots,N_m\in[\SS]$ and $\La^{(k)}\in N_k$ are arbitrary, then $(B,\La^{(1)}_n,\ldots,\La^{(m)}_n)\To (B,\II_{N_1}, \ldots,\II_{N_m})$ in $D_{\RR^{m+1}}[0,\infty)$.

It is easily verified that $[\SS]$ is a vector space under the operations $c[N]=[cN]$ and $[M]+[N]=[M+N]$, and that $N\mapsto \II_N$ is linear and injective. This gives us a one-to-one correspondence between $[\SS]$ and processes of the form \eqref{I_def}.

If $X=f(B)$, where $f\in C^\infty(\RR)$, then we define $N^X=[\La^X] \in[\SS]$. We may now identify $X$ with $N^X$, and will sometimes abuse notation by writing $X\in[\SS]$. It may therefore be necessary to deduce from context whether $X$ refers to the process $f(B)$, the sequence $\La^X=\{\La^X_n\}$, or the equivalence class $N^X = [\La^X]$. Typically, there will be only one sensible interpretation, but when ambiguity is possible, we will be specific.

Note that, using \eqref{smooth_expan}, we obtain $\ph_{1,X}=f'$, $\ph_{3,X}=\frac1{24}f'''$, $\II_X(0)=X(0)=f(0)$, and $\Phi_X=f-f(0)$. Hence, by \eqref{I_def}, we have $\II_X(t) = X(t)$. Because of this, and because of the one-to-one correspondence between $N\in[\SS]$ and the process $\II_N(t)$ in \eqref{I_def}, we will sometimes abuse notation and write $N(t)=N_t=\II_N(t)$. Again, when there is a possible ambiguity as to whether $N$ refers to an element of $[\SS]$ or to the process $\II_N$, we will be specific.

\subsection{The signed cubic variation}

If $\La\in\SS$, we define $V_n^\La(t)=\sum_{j=1}^{\flr{nt}}(\de_j(\La_n))^3$ and $V^\La=\{V^\La_n\}$. Since $\de_j(V^\La_n) = (\de_j(\La_n))^3$, it is easy to see from \eqref{script_S} that $V^\La\in\SS$, $\ph_{1,V^\La}=0$, and $\ph_{3, V^\La}=\ph_{1,\La}^3$. Hence, if $\La\equiv\Th$, then $V^\La\equiv V^\Th$. We may therefore define the signed cubic variation of $N =[\La]\in[\SS]$ to be $[V^\La]\in[\SS]$. We denote the signed cubic variation of $N$ by $\cub{N}$. We then have $\ph_{1,\cub{N}}=0$, $\ph_{3, \cub{N}}=\ph_{1,N}^3$, and $\II_{\cub{N}}(0)=0$, so that by \eqref{I_def},
  \[
  \cub{N}_t = \II_{\cub{N}}(t)
    = \ka\int_0^t (\ph_{1,N}(B(s)))^3\,dW(s).
  \]
For example, suppose $X=f(B)$, where $f\in C^\infty(\RR)$. Then $\cub{X} = \cub{N^X}$. Since $N^X=[\La^X]$, we have $\cub{N^X}=[V^{\La^X}]$. Note that $V^{\La^X}=\{V^{\La^X}_n\}$ and
  \[
  V^{\La^X}_n(t) = \sum_{j=1}^{\flr{nt}}(\de_j(\La^X_n))^3
    = \sum_{j=1}^{\flr{nt}}(X(t_j) - X(t_{j-1}))^3.
  \]
In other words, $\cub{X}$ is the equivalence class in $\SS$ of the above sequence of sums of cubes of increments of $X$. By Theorem \ref{T:gen_main}, $\cub{X}_t=\II_{\cub{X}}(t)$ is the stochastic process which is the limit in law of this sequence. Since $\ph_{1,X}=f'$, we have $\ph_{1,\cub{X}}=0$ and $\ph_{3,\cub{X}}=(f')^3$, so that
  \[
  \cub{X}_t = \II_{\cub{X}}(t) = \ka\int_0^t (f'(B(s)))^3\,dW(s).
  \]
In particular, taking $f(x)=x$ gives $\cub{B}_t=\ka W$.

\subsection{The weak Stratonovich integral}

If $\La_n,\Th_n\in\SS_n$, then we define
  \[
  (\La_n\circ\Th_n)(t) = \sum_{j=1}^{\flr{nt}}
    \frac{\La_n(t_{j-1}) + \La_n(t_j)}2\de_j(\Th_n).
  \]
If $\La,\Th\in\SS$, then we define $\La\circ\Th=\{\La_n\circ\Th_n\}_{n = 1}^\infty$.

\begin{lemma}\label{L:WSIC}
If $X=f(B)$, where $f\in C^\infty(\RR)$ and $\La\in\SS$, then $\La^X \circ\La\in\SS$. Moreover, if $\La\equiv\Th$, then $\La^X\circ\La \equiv\La^X\circ\Th$.
\end{lemma}

\pf Clearly, $\La^X_n\circ\La_n\in\SS_n$ and $\La^X_n\circ\La_n(0) =0$ for all $n$, so that Definition \ref{D:script_S}(i) and (ii) hold. For $a,b\in\RR$, we use the Taylor expansion
  \[
  \frac{f(b) + f(a)}2 = f(x) + \frac18f''(x)(b - a)^2
    + \frac1{4!2^4}f^{(4)}(x)(b - a)^4
    + h(a,b)(b - a)^6,
  \]
where $x=(a+b)/2$, and $|h(a,b)|\le M(a,b)=\sup_{x\in[a\wedge b,a\vee b]} |g^{(6)}(x)|$. For a derivation of this Taylor expansion, see the proof of Lemma 5.2 in \cite{NRS}.

Taking $a=B(t_{j-1})$ and $b=B(t_j)$ gives
  \begin{align*}
  \frac{\La_n^X(t_{j-1}) + \La_n^X(t_j)}2
    &= \frac{f(B(t_{j-1})) + f(B(t_j))}2\\
  &= f(\be_j) + \frac18f''(\be_j)\De B_{j,n}^2
    + \frac1{4!2^4}f^{(4)}(\be_j)\De B_{j,n}^4
    + R_{j,n},
  \end{align*}
where for each $T,K>0$, there exists a finite constant $C_{T,K}$ such that
  \[
  |R_{j,n}|1_{\{B^*(T)\le K\}} \le C_{T,K}|\De B_{j,n}|^6,
  \]
whenever $j/n\le T$. Choose $\ph_5\in C^\infty(\RR)$ such that
  \[
  \de_j(\La_n) = \ph_{1,\La}(\be_j)\De B_{j,n}
    + \ph_{3,\La}(\be_j)\De B_{j,n}^3
    + \ph_5(\be_j)\De B_{j,n}^5 + \wt R_{j,n},
  \]
where for each $T,K>0$, there exists a finite constant $\wt C_{T,K}$ such that
  \[
  |\wt R_{j,n}|1_{\{B^*(T)\le K\}} \le \wt C_{T,K}|\De B_{j,n}|^7,
  \]
whenever $j/n\le T$. Then
  \begin{align*}
  \de_j(\La^X_n\circ\La_n)
    &= \frac{\La_n^X(t_{j-1}) + \La_n^X(t_j)}2\de_j(\La_n)\\
  &= (f\ph_{1,\La})(\be_j)\De B_{j,n}
    + \left({\frac18f''\ph_{1,\La} + f\ph_{3,\La}}\right)
    (\be_j)\De B_{j,n}^3
    + h(\be_j)\De B_{j,n}^5 + \wh R_{j,n},
  \end{align*}
for an appropriately chosen smooth function $h$, and with $\wh R_{j,n}$ satisfying Definition \ref{D:script_S}(iii).

It follows that $\La^X\circ\La\in\SS$, and that $\ph_{1,\La^X\circ\La} = f\ph_{1,\La}$ and $\ph_{3,\La^X\circ\La} = \frac18f''\ph_{1,\La} + f\ph_{3,\La}$. This implies that if $\La\equiv\Th$, then $\La^X \circ\La \equiv\La^X\circ\Th$. \qed

\bigskip

If $X=f(B)$, where $f\in C^\infty(\RR)$, and $N=[\La]\in[\SS]$, we may now define $X\circ N=[\La^X\circ\La]$. Note that if $Y=g(B)$, where $g\in C^\infty$, and $M\in[\SS]$, then $(X+Y)\circ N=X\circ N + Y\circ N$ and $X\circ(N+M)=X\circ N+X\circ M$. From the proof of Lemma \ref{L:WSIC}, we have
  \begin{align}
  \II_{X\circ N}(0) &= 0,\label{WSIC1}\\
  \ph_{1,X\circ N} &= f\ph_{1,N},\label{WSIC2}\\
  \ph_{3,X\circ N} &= \frac18f''\ph_{1,N} + f\ph_{3,N}.\label{WSIC3}
  \end{align}
We may use these formulas, together with \eqref{I_def}, to calculate $\II_{X\circ N}$, given $f$, $\ph_{1,N}$, and $\ph_{3,N}$.

We now adopt some more traditional notation. If $X=f(B)$, where $f \in C^\infty$, and $N\in[\SS]$, then
  \[
  \int X\,\dd N = X \circ N \in [\SS],
  \]
and
  \[
  \int_0^t X(s)\,\dd N(s) = (X \circ N)_t = \II_{X\circ N}(t).
  \]
As we noted earlier, there is a one-to-one correspondence between $[\SS]$ and processes of the form \eqref{I_def}. We may therefore go back and forth between the above two objects according to what is more convenient at the time. We will use the shorthand notation $\dd M=X\,\dd N$ to denote the equality $M=\int X\,\dd N$.

Before investigating our change-of-variable formulas, let us first consider some examples.

\begin{expl}
Let $X=f(B)$ and $Y=g(B)$, where $f,g\in C^\infty(\RR)$. Then
  \[
  \int X\,\dd Y = X \circ Y = X \circ N^Y = [\La^X \circ \La^Y],
  \]
and $\La^X\circ\La^Y=\{\La_n^X\circ\La_n^Y\}$, where
  \begin{align*}
  (\La_n^X\circ\La_n^Y)(t) &= \sum_{j=1}^{\flr{nt}}
    \frac{\La^X_n(t_{j-1}) + \La^X_n(t_j)}2\de_j(\La^Y_n)\\
  &= \sum_{j=1}^{\flr{nt}}
    \frac{X(t_{j-1}) + X(t_j)}2(Y(t_j) - Y(t_{j-1})).
  \end{align*}
In other words, $\int X\,\dd Y$ is the equivalence class in $\SS$ of the above sequence of Stratonovich-type symmetric Riemann sums. Also, $\int_0^t X(s)\,\dd Y(s)=\II_{X\circ Y}(t)$, so that by Theorem \ref{T:gen_main}, $\int_0^t X(s)\,\dd Y(s)$ is the stochastic process which is the limit in law of this sequence.
\end{expl}

\begin{expl}\label{E:cub_form}
Again let $X=f(B)$ and $Y=g(B)$, where $f,g\in C^\infty(\RR)$. Then
  \[
  \int X\,\dd\cub{Y} = X \circ \cub{Y}
    = [\La^X \circ V^{\La^Y}]
  \]
and $\La^X\circ V^{\La^Y}=\{\La_n^X\circ V_n^{\La^Y}\}$, where
  \begin{align*}
  (\La_n^X\circ V_n^{\La^Y})(t) &= \sum_{j=1}^{\flr{nt}}
    \frac{\La^X_n(t_{j-1}) + \La^X_n(t_j)}2\de_j(V_n^{\La^Y})\\
  &= \sum_{j=1}^{\flr{nt}}
    \frac{X(t_{j-1}) + X(t_j)}2(Y(t_j) - Y(t_{j-1}))^3.
  \end{align*}
In other words, $\int X\,\dd\cub{Y}$ is the equivalence class in $\SS$ of the above sequence of sums, and $\int_0^t X(s)\,\dd\cub{Y}_s =\II_{X\circ\cub{Y}}(t)$ is the the limit in law of this sequence. Recall that $\ph_{1,\cub{Y}}=0$ and $\ph_{3,\cub{Y}}=(g')^3$. Hence, by \eqref{WSIC2} and \eqref{WSIC3}, we have $\ph_{1,X\circ\cub{Y}}=f\ph_{1,\cub{Y}}=0$ and $\ph_{3,X\circ\cub{Y}} =\frac18f''\ph_{1,\cub{Y}} + f\ph_{3,\cub{Y}}=f(g')^3$, so that by \eqref{I_def}, we have
  \begin{equation}\label{cub_form}
  \int_0^t X(s)\,\dd\cub{Y}_s
    = \ka\int_0^t f(B(s))(g'(B(s)))^3\,dW(s).
  \end{equation}
\end{expl}

\begin{expl}
For one last example, let $X=f(B)$, $Y=g(B)$, and $Z=h(B)$, where $f,g,h\in C^\infty(\RR)$, and let $N=\int Y\,\dd Z$. Then
  \[
  \int X\,\dd N = X \circ N = X \circ [\La^Y \circ \La^Z]
    = [\La^X \circ (\La^Y \circ \La^Z)],
  \]
and
  \begin{align*}
  (\La^X_n \circ (\La^Y \circ \La^Z)_n)(t)
    &= (\La^X_n \circ (\La^Y_n \circ \La^Z_n))(t)\\
  &= \sum_{j=1}^{\flr{nt}} \frac{X(t_{j-1}) + X(t_j)}2
    \frac{Y(t_{j-1}) + Y(t_j)}2(Z(t_j) - Z(t_{j-1})).
  \end{align*}
Hence, $\int X\,\dd N$ is the equivalence class in $\SS$ of the above sequence of sums, and $\int_0^t X(s)\,\dd N(s)$ is the limit in law of this sequence.
\end{expl}

\section{Change-of-variable formulas}\label{S:formulas}

We have already identified smooth functions of $B$ with their corresponding sequences in $\SS$, as well as with their equivalence classes in $[\SS]$. In this section, it will be helpful to do the same for $\FF_\infty^B$-measurable random variables, which can serve as initial values for the stochastic processes we are considering.

Let $\eta$ be an $\FF_\infty^B$-measurable random variable, let $\La^\eta_n(t)=\eta$ for all $t\ge0$, and let $\La^\eta=\{\La^\eta_n\}$. Since $\de_j( \La^\eta_n)=0$ for all $j$ and $n$, we have that $\La^\eta\in\SS$. We may therefore identify $\eta$ with $\La^\eta\in\SS$, and also with $N^\eta=[\La^\eta]\in[\SS]$. Note, then, that $\ph_{1,\eta}=\ph_{3,\eta} =0$, and $\eta(t)=N^\eta(t)=\II_{N^\eta}(t)=\eta$ for all $t\ge0$. Note also that $\int X\,\dd\eta=0$.

We begin with the following result, which tells us that every element of $[\SS]$ has a unique decomposition into the sum of a smooth function of $B$ and an integral against $\cub{B}$.

\begin{lemma}\label{L:decomp}
Each $N\in[\SS]$ can be written as $N=\eta+Y+V$, where $\eta$ is an $\FF_\infty^B$-measurable random variable, $Y=g(B)$ for some $g\in C^\infty(\RR)$, and $V=\int \th(B)\,\dd\cub{B}$ for some $\th\in C^\infty(\RR)$.

Suppose $N=\wt\eta+\wt Y + \wt V$ is another such representation, with $\wt Y=\wt g(B)$ and $\wt V=\int\wt\th(B)\,\dd\cub{B}$. Let $c=g(0)-\wt g(0)$. Then $\wt\eta=\eta+c$, $\wt g=g-c$, and $\wt\th=\th$. In particular, there is a unique such representation with $g(0)=0$.

An explicit representation is given by $\eta=N(0)=\II_N(0)$, $\th = \ph_{3,N} - \frac1{24}\ph_{1,N}''$ and $g$ chosen so that $g'=\ph_{1,N}$ and $g(0)=0$.
\end{lemma}

\pf Let $N\in[\SS]$. Let $\eta=N(0)$, $\th = \ph_{3,N} - \frac1{24} \ph_{1,N}''$ and choose $g$ so that $g' = \ph_{1,N}$ and $g(0) =0$. Let $Y=g(B)$ and $V=\int \th(B)\,\dd\cub{B}$. To prove that $N=\eta+ Y+V$, it will suffice to show that
  \begin{align*}
  N(t) &= \eta(t) + Y(t) + V(t)\\
  &= N(0) + g(B(t)) + \int_0^t \th(B(s))\,\dd\cub{B}_s.
  \end{align*}
But this follows immediately from \eqref{I_def} and \eqref{cub_form}.

Now suppose $N(t)=\wt\eta+\wt g(B(t)) +\int_0^t\wt\th(B(s))\,\dd \cub{B}_s$. Then $E[N(t)\mid\FF_\infty^B]=\eta+g(B(t))=\wt\eta + \wt g(B(t))$ a.s., which gives $\eta-\wt\eta+(g-\wt g)(B(t))=0$ a.s. for all $t\ge0$. Hence, there exists a constant $c\in\RR$ such that $g-\wt g=c$, and it follows that $\wt\eta = \eta+c$. We then have $\MM(t)=\int_0^t(\th-\wt\th)(B(s))\,dW(s)=0$ a.s., so that $E[\MM(t)^2 \mid\FF_\infty^B]=\int_0^t|(\th-\wt\th)(B(s))|^2\,ds=0$ a.s. for all $t\ge0$, which implies $\th=\wt\th$. \qed

\bigskip

We next verify that processes of the form $V=\int\th(B)\,\dd\cub{B}$ behave as we would expect them to in regards to integration.

\begin{lemma}\label{L:usual}
Let $X=f(B)$, where $f\in C^\infty(\RR)$, and let $\th\in C^\infty(\RR)$. If $\dd V=\th(B)\, \dd\cub{B}$, then $X\,\dd V=X\th(B)\, \dd\cub{B}$.
\end{lemma}

\pf Let $V=\int \th(B)\, \dd\cub{B}$, $U=\int X\th(B)\, \dd\cub{B}$, and $N=\int X\,\dd V$. Since $N(0)=U(0)=0$, it will suffice to show that $\ph_{1,U}=\ph_{1,N}$ and $\ph_{3,U}=\ph_{3,N}$. By Example \ref{E:cub_form}, $\ph_{1,V}=\ph_{1,U}=0$, $\ph_{3,V} = \th$, and $\ph_{3,U}=f\th$. On the other hand, by \eqref{WSIC2} and \eqref{WSIC3}, we have $\ph_{1,N}=f\ph_{1,V}=0$ and $\ph_{3,N}=\frac18f''\ph_{1,V} +f\ph_{3,V}=f\th$. \qed

\bigskip

We finally present our main result for doing calculations with the weak Stratonovich integral.

\begin{thm}\label{T:main}
Let $N\in[\SS]$ and write $N=\eta+Y+V$, where $\eta$ is an $\FF_\infty^B$-measurable random variable, $Y=g(B)$, and $V=\int\th(B)\,\dd\cub{B}$ for some $g,\th\in C^\infty(\RR)$. Let $X=f(B)$, where $f\in C^\infty(\RR)$. Then
  \begin{equation}\label{main_big}
  \int X\,\dd N = \Phi(B)
    + \frac1{12}\int(f''g' - f'g'')(B)\,\dd\cub{B}
    + \int X\,\dd V,
  \end{equation}
where $\Phi\in C^\infty(\RR)$ is chosen so that $\Phi' =fg'$ and $\Phi(0) =0$.
\end{thm}

\begin{rmk}
Since $M=\int X\,\dd N\in[\SS]$, Lemma \ref{L:decomp} tells us that $M$ has a unique decomposition into the sum of a smooth function of $B$ and an integral against $\cub{B}$. Theorem \ref{T:main} gives us a convenient formula for this decomposition.
\end{rmk}

\begin{rmk}
Theorem \ref{T:main}, and the corollaries that are to follow, express equalities in the space $[\SS]$. Each side of \eqref{main_big} is an equivalence class of sequences of Riemann sums that converge in law. The equivalence relation is such that if we choose any sequence from the class on the left and any sequence from the class on the right, then their difference will converge to zero ucp. Note that this is a stronger statement than simply asserting that the two sequences have the same limiting law.
\end{rmk}

\noindent\textbf{Proof of Theorem \ref{T:main}.} Since $\int X\,\dd N=\int X\,\dd\eta+\int X\,\dd Y+\int X\,\dd V$ and $\int X\,\dd\eta=0$, it follows from \eqref{cub_form} that we need only show
  \begin{equation}\label{main}
  \int_0^t X(s)\,\dd Y(s) = \Phi(B(t))
    + \frac\ka{12}\int_0^t (f''g' - f'g'')(B(s))\,dW(s).
  \end{equation}
By \eqref{I_def}, we have
  \[
  \int_0^t X(s)\,\dd Y(s) = \Phi_M(B(t))
    + \ka\int_0^t\bigg(\ph_{3,M} - \frac1{24}\ph_{1,M}''
    \bigg)(B(s))\,dW(s),
  \]
where $M = X\circ Y$. Recall that $\ph_{1,Y} =g'$ and $\ph_{3,Y} = \frac 1 {24} g'''$. By \eqref{WSIC2} and \eqref{WSIC3}, we have $\ph_{1,M}=fg'$ and $\ph_{3,M}=\frac18f''g'+\frac1{24}fg'''$. Since $\Phi_M(0)=0$ and $\Phi_M'=\ph_{1,M}=fg'$, we have $\Phi_M=\Phi$, and we also have
  \begin{align*}
  \ph_{3,M} - \frac1{24}\ph_{1,M}''
    &= \frac18 f''g' + \frac1{24} fg''' - \frac1{24}(fg')''\\
  &= \frac18 f''g' + \frac1{24} fg'''
    - \frac1{24}f''g' - \frac1{12}f'g'' - \frac1{24}fg'''\\
  &= \frac1{12}(f''g' - f'g''),
  \end{align*}
and this verifies \eqref{main}. \qed

\begin{cor}\label{C:ito_form_Y}
Let $Y=g(B)$, where $g\in C^\infty(\RR)$, and let $\ph\in C^\infty$. Then
  \begin{equation}\label{ito_form_Y}
  \ph(Y(t)) = \ph(Y(0)) + \int_0^t \ph'(Y(s))\,\dd Y(s)
    - \frac1{12}\int_0^t \ph'''(Y(s))\,\dd\cub{Y}_s.
  \end{equation}
\end{cor}

\pf Let $X=\ph'(Y)=f(B)$, where $f=\ph'\circ g$. By Theorem \ref{T:main},
  \[
  \int X\,\dd Y = \Phi(B)
    + \frac1{12}\int(f''g' - f'g'')(B)\,\dd\cub{B},
  \]
where $\Phi\in C^\infty(\RR)$ is chosen so that $\Phi' =fg'$ and $\Phi(0) =0$. Since $(\ph\circ g)'=fg'$, we have $\Phi=(\ph\circ g)-(\ph\circ g)(0)$. Also,
  \[
  f''g' - f'g'' = ((\ph'''\circ g)(g')^2 + (\ph''\circ g)g'')g'
    - (\ph''\circ g)g'g'' = (\ph'''\circ g)(g')^3.
  \]
Thus,
  \begin{align*}
  \int_0^t \ph'(Y(s))\,\dd Y(s) &= \int_0^t X(s)\,\dd Y(s)\\
  &= (\ph\circ g)(B(t)) - (\ph\circ g)(0)
    + \frac1{12}\int_0^t (\ph'''\circ g)(B(s))
    (g'(B(s)))^3\,\dd\cub{B}_s\\
  &= \ph(Y(t)) - \ph(Y(0))
    + \frac\ka{12}\int_0^t \ph'''(Y(s))(g'(B(s)))^3\,dW(s).
  \end{align*}
By \eqref{cub_form}, this gives
  \[
  \int_0^t \ph'(Y(s))\,\dd Y(s)
    = \ph(Y(t)) - \ph(Y(0))
    + \frac1{12}\int_0^t \ph'''(Y(s))\,\dd\cub{Y}_s,
  \]
which is \eqref{ito_form_Y}. \qed

\begin{cor}\label{C:sub_weak_semi}
Let $N\in[\SS]$ and write $N=\eta+Y+V$, where $\eta$ is an $\FF_\infty^B$-measurable random variable, $Y=g(B)$, and $V=\int\th(B)\,\dd\cub{B}$ for some $g,\th\in C^\infty(\RR)$. Let $X=f(B)$ and $Z=h(B)$, where $f,h\in C^\infty(\RR)$. Then:
  \begin{equation}\label{sub_weak_semi}
  \text{If $\dd M=X\,\dd N$, then
    $Z\,\dd M = ZX\,\dd N-\frac14(f'g'h')(B)\,\dd\cub{B}$.}
  \end{equation}
Moreover, the above correction term is a ``weak triple covariation" in the following sense: If $\VV=\{\VV_n\}$, where
  \[
  \VV_n(t) = \sum_{j=1}^{\flr{nt}}
    (X(t_j) - X(t_{j-1}))(Y(t_j) - Y(t_{j-1}))(Z(t_j) - Z(t_{j-1})),
  \]
then $\VV\in\SS$ and $[\VV]=\int(f'g'h')(B)\,\dd\cub{B}$.
\end{cor}

\pf Let $N$, $X$, and $Z$ be as in the hypotheses, and let $M=\int X\,\dd N$. By Theorem \ref{T:main},
  \[
  M = \Phi(B) + \frac1{12}\int(f''g' - f'g'')(B)\,\dd\cub{B}
    + \int X\,\dd V,
  \]
where $\Phi\in C^\infty(\RR)$ is chosen so that $\Phi'=fg'$ and $\Phi(0) = 0$. Hence, by Lemma \ref{L:usual},
  \begin{equation}\label{main_pf1}
  \int Z\,\dd M = \int Z\,\dd\Phi(B)
    + \frac1{12}\int(f''g'h - f'g''h)(B)\,\dd\cub{B}
    + \int ZX\,\dd V.
  \end{equation}
By Theorem \ref{T:main},
  \[
  \int Z\,\dd\Phi(B) = \Psi(B)
    + \frac1{12}\int(h''\Phi' - h'\Phi'')(B)\,\dd\cub{B},
  \]
where $\Psi\in C^\infty(\RR)$ is chosen so that $\Psi'=h\Phi'$ and $\Psi(0)=0$. Theorem \ref{T:main} also gives
  \[
  \int ZX\,\dd Y = \wt\Psi(B)
    + \frac1{12}\int((fh)''g' - (fh)'g'')(B)\,\dd\cub{B},
  \]
where $\wt\Psi\in C^\infty(\RR)$ is chosen so that $\wt\Psi'=fhg'$ and $\wt\Psi(0)=0$. Note, however, that this implies $\Psi =\wt\Psi$, which gives
  \[
  \int Z\,\dd\Phi(B) = \int ZX\,\dd Y + \frac1{12}\int
    (h''\Phi' - h'\Phi'' - (fh)''g' + (fh)'g'')(B)\,\dd\cub{B}.
  \]
Substituting $\Phi'=fg'$ into the above and simplifying gives
  \[
  \int Z\,\dd\Phi(B) = \int ZX\,\dd Y + \frac1{12}\int
    (f'g''h - f''g'h - 3f'g'h')(B)\,\dd\cub{B}.
  \]
Substituting this into \eqref{main_pf1} gives
  \begin{align*}
  \int Z\,\dd M &= \int ZX\,\dd Y
    - \frac14\int(f'g'h')(B)\,\dd\cub{B} + \int ZX\,\dd V\\
  &= \int ZX\,\dd N - \frac14\int(f'g'h')(B)\,\dd\cub{B},
  \end{align*}
and this verifies \eqref{sub_weak_semi}.

Finally, if $\VV=\{\VV_n\}$, then $\de_j(\VV_n) = \de_j(\La^X_n) \de_j(\La^Y_n)\de_j(\La^Z_n)$. From \eqref{script_S}, we see that $\VV\in\SS$, $\ph_{1,\VV}=0$, and $\ph_{3,\VV}=\ph_{1,X}\ph_{1,Y} \ph_{1,Z}=f'g'h'$. Since $[\VV]_0=0$, it follows from Example \ref{E:cub_form} that $[\VV]=\int(f'g'h')(B)\,\dd\cub{B}$. \qed

\section*{Acknowledgments}
\addcontentsline{toc}{section}{Acknowledgments}

Thanks go to Tom Kurtz and Frederi Viens for stimulating and helpful comments, feedback, and discussions.


%



\end{document}